\documentclass[11pt]{article}
\usepackage{amsmath,amsthm,amsfonts,verbatim}

\newcommand{\R}{\mathbb{R}}

\theoremstyle{plain}
 \newtheorem{theorem}{Theorem}
 \newtheorem{lemma}{Lemma}
 \newtheorem{corollary}[lemma]{Corollary}
 \newtheorem{remark}{Remark}
 
 \newtheorem{proposition}[lemma]{Proposition}

 \theoremstyle{definition}
 \newtheorem{definition}{Definition}

\begin{document}

\title{On Some Low Distortion Metric Ramsey Problems}

\date{}

\author{Yair Bartal\thanks{Supported in part by a grant from the
Israeli National Science Foundation.} \and Nathan Linial\thanks{Supported in
part by a grant from the Israeli National Science Foundation.} \and Manor
Mendel\thanks{Supported in part by the Landau Center.} \and Assaf Naor }

\maketitle

\begin{abstract}
In this note, we consider the metric Ramsey problem for the normed
spaces $\ell_p$. Namely, given some $1\le p \le \infty$ and
$\alpha \ge 1$, and an integer $n$, we ask for the largest $m$
such that every $n$-point metric space contains an $m$-point
subspace which embeds into $\ell_p$ with distortion at most $
\alpha$. In \cite{blmn0} it is shown that in the case of $\ell_2$,
the dependence of $m$ on $\alpha$ undergoes a phase transition at
$\alpha =2$. Here we consider this problem for other $\ell_p$, and
specifically the occurrence of a phase transition for $p\neq 2$.
It is shown that a phase transition does occur at $\alpha=2$ for every $p\in
[1,2]$. For $p>2$ we are unable to determine the answer, but
estimates are provided for the possible location of such a phase
transition. We also study the analogous problem for isometric
embedding and show that for every $1<p<\infty$ there are
arbitrarily large metric spaces, no four points of which embed
isometrically in $\ell_p$ .
\end{abstract}

\section{Introduction}\label{introduction}

A Ramsey-type theorem states that large systems necessarily contain large,
highly structured sub-systems. Here we consider Ramsey-type problems for
finite metric spaces, interpreting  ``highly structured" as having low
distortion embedding in $\ell_p$.

A mapping between two metric spaces $f:M
\rightarrow X$, is called an embedding of $M$ in $X$.
The \emph{distortion} of the embedding is defined as
\[
\mathrm{dist}(f)=\sup_{\substack{x,y\in M\\x\neq
y}}\frac{d_X(f(x),f(y))}{d_M(x,y)}\cdot \sup_{\substack{x,y\in M\\x\neq
y}}\frac{d_M(x,y)}{d_X(f(x),f(y))}.
\]
The least distortion required to embed $M$ in $X$ is denoted by
$c_X(M)$. When $c_X(M)\leq \alpha$ we say that $M$ $\alpha$-embeds
in $X$. In  this note we study the following notion.

\begin{definition}[Metric Ramsey function]
We denote by
$R_{X}(\alpha,n)$ the largest integer $m$ such
that every $n$-point metric space has a subspace of size $m$ that
$\alpha$-embeds into $X$.
\end{definition}
 When $X=\ell_p$ we use the notations $c_p$ and $R_{p}$. Note that for $p\in [1,\infty]$, it is always true that
$R_p(\alpha,n)\geq R_2(\alpha,n)$. When $\alpha=1$ we drop it from the
notation, i.e., $R_X(n)=R_X(1,n)$.

Bourgain, Figiel, and Milman~\cite{bfm}  study this function for
$X=\ell_2$, as a metric analog of Dvoretzky's theorem~\cite{dvo}.
They prove

\begin{theorem}[\cite{bfm}] \label{thm:bfm86}
For any $\alpha>1$ there exists $C(\alpha)>0$ such that $R_2(\alpha,n) \geq
C(\alpha) \log n$. Furthermore, there exists $\alpha_0>1$ such that
$R_2(\alpha_0,n)=O(\log n)$.
\end{theorem}

In~\cite{blmn0} the metric Ramsey problem is studied
comprehensively. In particular, the following
phase transition is established in the case of $X=\ell_2$.

\begin{theorem}[\cite{blmn0}]
\label{thm:phase} Let  $n\in \mathbb{N}$. Then:
\begin{enumerate}
\item\label{item:upper} For every $1<\alpha<2$: \( c(\alpha)\log
n\le R_2(\alpha,n)\le 2\log n+C(\alpha), \) where
$c(\alpha),C(\alpha)$ may depend only on $\alpha$.

\item\label{item:lower} For every $\alpha>2$: \( n^{c'(\alpha)}\le
R_2(\alpha,n)\le n^{C'(\alpha)}, \) where $c'(\alpha),C'(\alpha)$
depend only on $\alpha$ and $0<c'(\alpha)\leq C'(\alpha)<1$.
Moreover, $c'(\alpha)$ tends to $1$ as $\alpha$ tends to $\infty$.
\end{enumerate}
\end{theorem}

By Dvoretzky's theorem, the lower bound in part~\ref{item:lower}
of Theorem~\ref{thm:phase} implies in particular that if
$\alpha>2$, and $X$ is any infinite dimensional normed space, then
$R_X(\alpha,n)\ge n^{c'(\alpha)}$. Therefore, in our search
for a possible phase transition for $R_p(\cdot,n)$, $p\neq
2$, it is natural to extend the upper bound in
part~\ref{item:upper} of Theorem~\ref{thm:phase} to this range.
The main result proved in this note is the following:

\begin{theorem}\label{thm:<2 a} There is an absolute constant $c>0$ such that
for every $0<\delta<1$,
\begin{enumerate}
\item For $1\le p< 2$, \( R_p({2-\delta},n)\le
e^{\frac{c}{\delta^2}}\log n. \) \item For $2<p<\infty$, $
R_p({2^{2/p}-\delta},n)\le e^{\frac{c}{p^2\delta^2}}\log n. $
\end{enumerate}
\end{theorem}

Thus we extend the result of \cite{blmn0} to show that a phase
transition occurs in the metric Ramsey problem for $\ell_p$,
$p\in[1,2)$, at $\alpha=2$. The asymptotic behavior of
$R_p(\alpha,n)$ for $p>2$, and $\alpha \in [2^{2/p},2]$, is left
as an open problem. In particular, we do not know whether or not
this function undergoes a similar phase transition. We find this
problem potentially significant: if there is a phase transition at
$2$ also in the range $2<p<\infty$, then this result will
certainly be of great interest. On the other hand, if it is
possible to improve the lower bound in part~\ref{item:lower} of
Theorem~\ref{thm:phase} for $p>2$ and certain distortions strictly
less than $2$, then this would involve an embedding technique that
is different from the method used in~\cite{blmn0}, which doesn't
distinguish between the various $\ell_p$ spaces.


The proof of the upper bound on $R_2(\alpha,n)$ for $\alpha<2$
stated in Theorem~\ref{thm:phase} uses the Johnson-Lindenstrauss
dimension reduction lemma for $\ell_2$ \cite{jl}. For $\ell_p$,
$p\neq 2$, no such dimension reduction is known to hold. (Recent
work~\cite{br:cha, lee:naor} shows that dimension reduction does
not, in genereal, hold in $\ell_1$.) Our proof is based on a
non-trivial modification of the random construction in~\cite{bfm},
in the spirit of Erd\"os' upper bound on the Ramsey
numbers~\cite{erdos, bollobasrandom}.  In the process we prove
tight bounds on the embeddability of the metrics of complete
bipartite graphs in $\ell_p$. Specifically we show that
\[ c_p(K_{n,n})=\begin{cases} 2-\Theta(n^{-1}) & p\in[1,2] \\
2^{2/p} -\Theta((pn)^{-1}) & p>2. \end{cases}
\]

\medskip
The second part of this note addresses the isometric Ramsey
problem for $p\in (1,\infty)$. It turns out that this problem is
naturally tackled within the class of uniformly convex normed
spaces (see Section~\ref{section:isometric} for the definition).

\begin{theorem}[Isometric Ramsey Problem]\label{thm:isometric-lp}
Let $X$ be a uniformly convex normed space with
$\mathrm{dim}(X)\ge 2$. Then $R_X(1,n) = 3$ for $n\ge 3$.
\end{theorem}
Since $\ell_p$ is uniformly convex for $p\in (1,\infty)$, the
conclusion of Theorem~\ref{thm:isometric-lp} holds in these cases.
Note that the theorem does not apply for $\ell_1$ and
$\ell_\infty$ which are not uniformly convex. Specifically, it is
known that $\ell_\infty$ is universal in that it contains an
isometric copy of every finite metric space, whence
$R_\infty(n)=n$. It is known \cite{dezalau} that any 4-point
metric space is isometrically embeddable in $\ell_1$, and
therefore $R_1(n)\geq 4$ for $n \ge 4$. The determination of
$R_1(n)$ is left as an open problem.

\section{An Upper Bound For $\alpha<2$} \label{section:a<2}

In this section we prove that for any $\alpha<\min\{2,2^{2/p}\}$,
$R_p(\alpha,n)=O(\log n)$. Our technique both improves and
simplifies the technique of \cite{bfm}, which is itself in the
spirit of Erd\"os' original upper bound for the Ramsey coloring
numbers. The basic idea is to exploit a
universality property of random graphs $G\in G(n,1/2)$.
Namely, that any fixed graph of constant size appears
as an induced subgraph of every induced subgraph of $G$ of size
$\Omega(\log n)$. More precisely, we define the following notion
of universality.

\begin{definition} Let $H$ be a graph. A graph $G$ is called
$({H},s)$-universal if every set of $s$ vertices in $G$ contains
an induced subgraph isomorphic to $H$.
\end{definition}

\begin{proposition} \label{prop:ramsey-graph}
For every $k$-vertex graph ${H}$ there exists a constant $C>0 $
and an integer $n_0$ such that for any $n>n_0$ there exists a
$({H},C\log n)$-universal graph on $n$ vertices. Furthermore,
\[
C\le O\left(k^2 2^{\binom{k}{2}}\right)\quad \mathrm{and} \quad n_0 \le
O\left(k^3 2^{\binom{k}{2}}\right).
\]
\end{proposition}

Such facts are well-known in random graph theory, and similar
arguments can be found for example in \cite{ly}. We sketch the
standard details for the sake of completeness.

Recall that a family of sets $\mathcal{F}$ is called \emph{almost
disjoint} if $|A \cap B| \leq 1$ for every $A,B\in \mathcal{F}$.
In what follows, given a set $S$ and an integer $k$, we denote by
$\binom{S}{k}$ the set of all $k$-point subsets of $S$.
\begin{lemma} \label{prop:pairwise}
For every integer $k$ and a finite set $S$ of cardinality
$s=|S| > 2k^2$,
there exists an almost disjoint family $K\subset \binom{S}{k}$,
such that $|K|\geq \left\lfloor \tfrac{s}{2k}\right\rfloor ^2$.
\end{lemma}
\begin{proof}
Let $p$ be a prime satisfying $\tfrac{s}{2k}\leq p\leq
\tfrac{s}{k}$, and assume that
 \[ L=\{(i,j); i,j \in \mathbb{Z}_p,\ i\in\{0,\ldots,k-1\} \}
 \subseteq S.\]
For each $a,b\in \mathbb{Z}_p$ (the field of residues modulo $p$), define
\[ A_{a,b}=\{(i,j); \ j\equiv ai+b \pmod{p},\ i\in\{0,\ldots,k-1\}\} ,\]
and take $K=\{A_{a,b} | a,b\in \mathbb{Z}_p\}$. The set $K$ is
easily checked to satisfy the requirements.
\end{proof}

As usual $G(n,1/2)$ denotes the model of random graphs in
which each edge on $n$ vertices is chosen independently with
probability $1/2$.

\begin{lemma} \label{lem:induced}
Let $H$ be a $k$-vertex graph and let $s>2k^2$. The probability that
a random graph $G \in G(s,1/2)$ does not contain an
induced subgraph isomorphic to $H$, is at most
$(1-2^{-\binom{k}{2}})^{\left\lfloor\tfrac{s}{2k}\right\rfloor^2} $.
\end{lemma}
\begin{proof}
Construct, as in Lemma~\ref{prop:pairwise}, an almost disjoint
family $\mathcal{F}$ of $\left\lfloor\tfrac{s}{2k}
\right\rfloor^2$ subsets of $\{1,\ldots,s\}$, the vertex set of
$G$. If $F_1 \neq F_2 \in \mathcal{F}$, then the event that the
restriction of $G$ to ${F_1}$ (resp. $F_2$) is isomorphic to $H$
are independent. Hence, the probability that none of the sets $F
\in \mathcal{F}$ spans a subgraph isomorphic to $H$ is at most
$(1-2^{-\binom{k}{2}})^{\left\lfloor\tfrac{s}{2k}\right\rfloor^2}
$.
\end{proof}

\begin{proof} [Proof of Proposition~\ref{prop:ramsey-graph}]
Let $G$ be a random graph in $G(n,1/2)$. By the previous lemma,
the expected number of sets of $s$ vertices which contain no
induced isomorphic copy of $H$ is at most $\binom{n}{s}
\left(1-2^{-\binom{k}{2}}\right)^{\left\lfloor\tfrac{s}{2k}\right\rfloor^2}$.
If this number is $< 1$, then there is an $(H,s)$-universal graph,
as claimed. It is an easy matter to check that this holds with the
parameters as stated.
\end{proof}

A class $\cal C$ of finite metric spaces is called a \emph{metric
class} if it is closed under isometries. $\cal C$ is said to be
\emph{hereditary}, if $M\in \cal C$ and $N\subset M$ imply $N\in
\cal C$. We call a metric space $(X,d)$ a $\{0,1,2\}$ metric space
if for all $x,y\in X$, $d(x,y)\in \{0,1,2\}$. There is a simple
1:1 correspondence between graphs and $\{0,1,2\}$ metrics. Namely,
associated with a $\{0,1,2\}$ metric space $M=(X,d)$ is the graph
$G=(X,E)$ where $\{x,y\} \in E$ iff $d_M(x,y)=1$.

\begin{lemma} \label{lem:p3}
Let $\cal C$ be a hereditary metric class of finite metric spaces,
and suppose that there exists some finite $\{0,1,2\}$ metric space
${M}_0$ which is not in $P$. Then there exist metric spaces
$M=M_n$ of arbitrarily large size $n$ such that every subspace
$S\subset M_n$ with at least $C \log n$ points is not in $\cal C$.
The constant $C$ depends only on the cardinality of $M_0$.
\end{lemma}
\begin{proof}
Let ${H}_0$ be the graph corresponding to the metric space
${M}_0$.
We apply
Proposition~\ref{prop:ramsey-graph}, to construct arbitrarily
large graphs $G_n=(V_n,E_n)$ with $|V_n|=n$, in which every
set of $\ge C\log n$ vertices contains an
induced subgraph isomorphic to $H_0$. Let $M_n$ be the $n$-point
metric space corresponding to $G_n$. It follows that
every subspace of $M_n$ of size $\ge C \log n$ contains
a metric subspace that is isometric to $M_0$.
Since $\cal C$ is hereditary, $S\notin \cal C$.
\end{proof}

Note that $\{M; M \text{ is a metric space, } c_p(M)\leq \alpha\}$
is a hereditary metric class. Therefore, in order to show that for
$\alpha<2$, $R_p(\alpha,n)=O(\log n)$, it is enough to find a
$\{0,1,2\}$ metric space whose $\ell_p$ distortion is greater than
$\alpha$. We use the complete bipartite graphs $K_{n,n}$. The
$\ell_p$-distortion of $K_{n,n}$, $1\le p<\infty$, is estimated in
the following proposition.

\begin{proposition}\label{prop:knn} For every $1\le p\le 2$,
$$
2\left(\frac{n-1}{n}\right)^{1/p}\le c_p(K_{n,n})\le
2\sqrt{\frac{n-1}{n}}
$$
For every $2\le p<\infty$,
$$
2^{2/p}\left(\frac{n-1}{n}\right)^{1/p}\le c_p(K_{n,n})\le
2^{2/p}\left(1-\frac{1}{2n}\right) ^{1/p}.
$$
\end{proposition}

Before proving Proposition~\ref{prop:knn}, we will deduce the main
result of this section:

\begin{theorem}\label{thm:<2} There is an absolute constant $c>0$ such that
for every $0<\delta<1$, if $1\le p\le 2$ then:
$$
R_p({2-\delta},n)\le e^{\frac{c}{\delta^2}}\log n,
$$
and if $2<p<\infty$ then:
$$
R_p({2^{2/p}-\delta},n)\le e^{\frac{c}{p^2\delta^2}}\log n.
$$
\end{theorem}
\begin{proof}

Proposition~\ref{prop:ramsey-graph} implies that
there is an absolute constant $C$ such that for every $n\ge
2^{Ck^3}$ there exists a $\{0,1,2\}$ metric space $M_n$ such that
any subset $S\subset M_n$ of cardinality at least $2^{Ck^2}\log n$
contains an isometric copy of $K_{k,k}$.

We start with $1\le p\le 2$. Let
$k=\left\lfloor\frac{2}{\delta}\right\rfloor+1$.
By Proposition~\ref{prop:knn},
$$
c_p(K_{k,k})\ge 2\left(1-\frac{1}{k}\right)^{1/p}>
2\left(1-\frac{\delta}{2}\right)=2-\delta,
$$
so that for $n$ large enough ($\ge e^{\frac{C'}{\delta^3}}$),
and hence for all $n$ (by proper choice of constants),
$$
R_p(2-\delta,n)\le e^{\frac{C'}{\delta^2}}\log n.
$$

When $p>2$ take $k=2\left\lfloor\frac{4}{p\delta}\right\rfloor$.
In this case one easily verifies that:
$$
c_p(K_{k,k})\ge 2^{2/p}\left(1-\frac{1}{k}\right)^{1/p}\ge
2^{2/p}-\delta,
$$
from which the required result follows as above.
\end{proof}

\medskip

In order to prove Proposition~\ref{prop:knn}, we need some
preparation.

\begin{lemma}\label{lem:matrix}
Let $A=(a_{ij})$ be an $n\times n$ matrix and $2\le p<\infty$.
Then:
$$
\sum_{i=1}^n\sum_{j=1}^n\left(\left|\sum_{k=1}^n
a_{ik}-\sum_{k=1}^n a_{jk}\right|^p+ \left|\sum_{k=1}^n
a_{ki}-\sum_{k=1}^n a_{kj}\right|^p\right)\le
 \frac{(2n)^p}{2}\sum_{i=1}^n\sum_{j=1}^n |a_{ij}|^p.
$$
\end{lemma}
\begin{proof} We identify $\ell_p^{n^2}$ with
the space of all $n\times n$ matrices
$A=(a_{ij})$, equipped with the $\ell_p$ norm:
$$
\|A\|_p=\left(\sum_{i=1}^n\sum_{j=1}^n |a_{ij}|^p\right)^{1/p}.
$$
Define a linear operator $T:\R^{n^2}\to \R^{n^2}\oplus\R^{n^2}$
by:
$$
T(a_{ij})=\left( \sum_{k=1}^n a_{ik}-\sum_{k=1}^n
a_{jk}\right)_{ij}\oplus \left(\sum_{k=1}^n a_{ki}-\sum_{k=1}^n
a_{kj}\right)_{ij}.
$$
For $q\ge 1$ denote $\|T\|_{q\to q}=\max_{A\neq
0}\|T(A)\|_q/\|A\|_q$. Our goal is to show that $\|T\|_{p\to p}\le
2^{1-1/p}n$. By a result from the complex interpolation theory for
linear operators (see \cite{bl}), for $2\le p\le \infty$,
$\|T\|_{p\to p}\le \|T\|_{2\to 2}^{2/p}\cdot \|T\|_{\infty\to
\infty}^{1-2/p}$. It is therefore enough to prove the required
estimate for $p=2$ and $p=\infty$. The case $p=\infty$ is simple:
$$
\|T(A)\|_{\infty}=\max_{1\le i,j\le n}\max\left\{ \left|
\sum_{k=1}^n a_{ik}-\sum_{k=1}^n a_{jk}\right|, \left|\sum_{k=1}^n
a_{ki}-\sum_{k=1}^n a_{kj}\right|\right\}\le 2n\|A\|_{\infty}.
$$
For $p=2$ we have to show that:
$$
\sum_{i=1}^n\sum_{j=1}^n\left(\left|\sum_{k=1}^n
a_{ik}-\sum_{k=1}^n a_{jk}\right|^2+ \left|\sum_{k=1}^n
a_{ki}-\sum_{k=1}^n a_{kj}\right|^2\right)\le
 2n^2\sum_{i=1}^n\sum_{j=1}^n |a_{ij}|^2.
$$
This inequality follows from the following elementary identity:
{\setlength\arraycolsep{1pt}
\begin{eqnarray*}
2n^2\sum_{i=1}^n\sum_{j=1}^na_{ij}^2
&=&\sum_{i=1}^n\sum_{j=1}^n\left[\left(\sum_{k=1}^na_{ik}-\sum_{k=1}^na_{jk}\right)^2+
\left(\sum_{k=1}^na_{ki}-\sum_{k=1}^na_{kj}\right)^2\right]+\\
&\phantom{\le}&+
2\sum_{i=1}^n\sum_{j=1}^n\left(na_{ij}-\sum_{k=1}^na_{ik}-\sum_{k=1}^na_{kj}\right)^2.
\end{eqnarray*}}
\end{proof}

\begin{corollary}\label{corr:roundness} Let $1\le p <\infty$ and $x_1,\ldots,x_n, y_1\ldots
y_n\in \ell_p$. Then if $2\le p<\infty$,
$$
\sum_{i=1}^n\sum_{j=1}^n
\left(\|x_i-x_j\|_p^p+\|y_i-y_j\|_p^p\right)\le
2^{p-1}\sum_{i=1}^n\sum_{j=1}^n \|x_i-y_j\|_p^p.
$$
If $1\le p\le 2$ then:
$$
\sum_{i=1}^n\sum_{j=1}^n
\left(\|x_i-x_j\|_p^p+\|y_i-y_j\|_p^p\right)\le
2\sum_{i=1}^n\sum_{j=1}^n \|x_i-y_j\|_p^p.
$$
\end{corollary}
\begin{proof} By summation it is clearly enough to prove these inequalities
for $x_1,\ldots,x_n,y_1,\ldots,y_n\in \R$. If $2\le p<\infty$ then
the required result follows from an application of
Lemma~\ref{lem:matrix} to the matrix $a_{ij}=x_i-y_j$. If $1\le
p\le 2$ then consider $\ell_p$ equipped with the metric
$d(x,y)=\|x-y\|_p^{p/2}$. It is well known (see \cite{ww}) that
$(\ell_p,d)$ embeds isometrically in $\ell_2$, so that the case
$1\le p \le 2$ follows from the case $p=2$.

\end{proof}

\medskip\noindent{\bf Remark.}
In \cite{enfloroundness} P. Enflo defined the notion on generalized roundness of a metric
space. A metric space $(M,d)$ is said to have generalized
roundness $q\ge 0$ if for every $x_1,\ldots,x_n,y_1,\ldots, y_n\in
M$,
$$
\sum_{i=1}^n\sum_{j=1}^n (d(x_i,x_j)^q+d(y_i,y_j)^q)\le
2\sum_{i=1}^n\sum_{j=1}^nd(x_i,y_j)^q.
$$
Enflo proved that Hilbert space has generalized roundness $2$ and
in \cite{ltw} the concept of generalized roundness was
investigated and was shown to be equivalent to the notion of
negative type (see \cite{dezalau,ww} for the definition).
Particularly, it was proved in \cite{ltw} that for $1\le p<2$,
$\ell_p$ has generalized roundness $p$, which is precisely the
second statement in Corollary \ref{corr:roundness}. For the case $p=1$
simpler, more direct proofs can be given
which do not use reduction to the case $p=2$,
see e.g.~\cite{dezalau}. Observe that
Lemma~\ref{lem:matrix} would follow simply by convexity had it not
been for the additional factor $1/2$ on the right-hand side. This
factor is crucial for our purposes, and this is why the
interpolation argument was needed.

\begin{proof}[Proof of Proposition~\ref{prop:knn}.] We identify $K_{n,n}$
with the metric on $\{u_1,\ldots,u_n,v_1,\ldots,v_n\}$ where
$d(u_i,u_j)=d(v_i,v_j)=2$ for all $i\neq j$, and
$d(u_i,v_j)=1$ for every $1\le i,j\le n$.
Fix some $1\le p<\infty$ and let
$f:\{u_1,\ldots,u_n,v_1,\ldots,v_n\}\to \ell_p$ be an embedding
such that for every $x,y\in K_{n,n}$, $d(x,y)\le
\|f(x)-f(y)\|_p\le Ld(x,y)$. Then,
$$
\sum_{i=1}^n\sum_{j=1}^n(\|f(u_i)-f(u_j)\|_p^p+\|f(v_i)-f(v_j)\|_p^p)\ge
2n(n-1)2^p,
$$
and
$$
\sum_{i=1}^n\sum_{j=1}^n\|f(u_i)-f(v_j)\|_p^p\le n^2L^p.
$$
For $1\le p\le 2$ Corollary \ref{corr:roundness} gives:
$$
2n(n-1)^p2^p\le 2n^2L^p\Longrightarrow L\ge
2\left(\frac{n-1}{n}\right)^{1/p}.
$$
For $2\le p<\infty$ we get that:
$$
2n(n-1)2^p\le 2^{p-1}n^2L^p\Longrightarrow L\ge
2^{2/p}\left(\frac{n-1}{n}\right)^{1/p}.
$$
This proves the required lower bounds on $c_p(K_{n,n})$.

To prove the upper bound assume first that $p=2$ and denote by
$\{e_i\}_{i=1}^{\infty}$ the standard unit vectors in $\ell_2$.
Define $f:K_{n,n}\to \ell_2^{2n}$ by:

\begin{align*}
f(u_i) &=\sqrt{2}\Bigl(e_i-\frac{1}{n}\sum_{j=1}^ne_j\Bigr),\\
f(v_i)
&=\sqrt{2}\Bigl(e_{n+i}-\frac{1}{n}\sum_{j=1}^ne_{n+j}\Bigr).
\end{align*}
Then for $i\neq j$,
$\|f(u_i)-f(u_j)\|_2=\|f(v_i)-f(v_j)\|_2=2=d(u_i,u_j)=d(v_i,v_j)$.
On the other hand: {\setlength\arraycolsep{1pt}
\begin{eqnarray*}
\|f(u_i)-f(v_j)\|_2&=&\sqrt{\|f(u_i)\|_2^2+\|f(v_j)\|_2^2}\\
&=&\sqrt{4\left(1-\frac{1}{n}\right)^2+4(n-1)\cdot\frac{1}{n^2}}=2\sqrt{\frac{n-1}{n}}.
\end{eqnarray*}}
This finishes the calculation of $c_2(K_{n,n})$. For $1\le p<2$,
since for every $\epsilon>0$ and for every $k$, $\ell_p$ contains
a $(1+\epsilon)$ distorted copy of $\ell_2^k$, we get the estimate
$c_p(K_{n,n})\le 2\sqrt{\frac{n-1}{n}}$.

The case $2<p< \infty$ requires a different embedding. We begin by
describing an embedding with distortion $2^{2/p}$ and then explain
how to modify it so as to reduce the distortion by a factor of
$\left(1-\frac{1}{2n}\right)^{1/p}$. Let $z_1,\ldots,z_n$ be a
collection of $n$ mutually orthogonal $\pm 1$ vectors of dimension
$m =O(n)$. (For example the first $n$ rows in an $m \times m$
Hadamard matrix). In our first embedding we define $f(u_i)$ as the
$(2m)$-dimensional vector $(z_i,0)$, namely, $z_i$ concatenated
with $m$ zeros. Likewise, $f(v_i)= (0,z_i)$ for all $i$. Now
$\|f(u_i) - f(u_j)\|_p = 2 \left(\frac{m}{2}\right)^{1/p}$ and
$\|f(u_i) - f(v_j)\|_p = (2m)^{1/p}$, and so $f$ has distortion
$2^{2/p}$. To get the $\left(1-\frac{1}{2n}\right)^{1/p}$
improvement, note that for some $m\le 4n$ it is possible to select
the $z_i$ so that the $m$-th coordinate in all of them is $+1$.
Modify the previous construction to an embedding into $2m-1$
dimensions as follows: Now $g(u_i)$ is $z_i$ concatenated with
$m-1$ zeros, whereas $g(v_i)$ has zeros in the first $m-1$
coordinates, $1$ in the $m$-th and this is followed by the first
$m-1$ coordinates of the vector $z_i$. The easy details are
omitted.
\end{proof}

\paragraph{Remark:} The upper bounds in Proposition~\ref{prop:knn} were not used in the proof of
Theorem~\ref{thm:<2}. Apart from their intrinsic interest, these
upper estimates show that the above technique cannot prove an
upper bound of $O(\log n)$ on $R_2({2-\epsilon},n)$ which is
independent of $\epsilon$. In fact, this can never be achieved
using $\{0,1,2\}$ metric spaces due to the following proposition.

\begin{proposition}
Let $X$ be an $n$-point $\{0,1,2\}$ metric space. Then $c_2(X)\le
2\sqrt{\frac{n-1}{n}}$.
\end{proposition}
\begin{proof}
We think of $X$ as a metric on $\{1,\ldots,n\}$ and denote
$d(i,j)=d_{ij}$. Define an $n\times n$ matrix $A=(a_{ij})$ as
follows:
 \[ a_{ij}=\begin{cases} 2 & \text{if}\ i=j \\
                                           0 & \text{if}\  d_{ij}=2\\
                                           \frac{2}{n} & \text{if}\ d_{ij}=1
                            \end{cases} .\]

We claim that $A$ is positive semidefinite. Indeed, for any $z\in
\R^n$ {\setlength\arraycolsep{1pt}
\begin{eqnarray*}
\langle Az,z\rangle&=&\sum_{i=1}^n\sum_{j=1}^n a_{ij}z_iz_j\\
&\ge& \sum_{i=1}^n 2z_i^2-\sum_{i \neq j}\frac{2}{n}|z_i|\cdot|z_j|\\
&\ge&  \sum_{i=1}^n 2z_i^2-\sum_{i=1}^n\sum_{j=1}^n\frac{2}{n}|z_i|\cdot|z_j|\\
&=&2\|z\|_2^2-\frac{2}{n}\|z\|_1^2\ge
2\|z\|_2^2-\frac{2}{n}n\|z\|_2^2=0.
\end{eqnarray*}}
In particular it follows that $A$ has a square root, denoted
$A^{1/2}$. Let $e_1,\ldots,e_n$ be the standard unit vectors in
$\R^n$. Define $f:X\to \R^n$ by $f(i)=A^{1/2}e_i$. Now,
$$
\|f(i)-f(j)\|_2^2=\langle Ae_i,e_i\rangle+\langle
Ae_j,e_j\rangle-2\langle Ae_i,e_j\rangle=a_{ii}+a_{jj}-2a_{ij},
$$
so that if $d_{ij}=1$ then $\|f(i)-f(j)\|_2=\sqrt{4-\frac{4}{n}}$
and if $d_{ij}=2$ then $\|f(i)-f(j)\|_2=2$. It follows that
$$
\mathrm{dist}(f)=2\sqrt{\frac{n-1}{n}}.
$$
\end{proof}

\section{The Isometric Ramsey Problem}\label{section:isometric}

In this section we prove that for $n\ge 3$, $1<p<\infty$,
$R_p(n)=R_p(1,n)=3$. In
fact, we show that this is true for any uniformly convex normed space. We
begin by sketching an argument that is specific to $\ell_2$:

\begin{proposition}
$R_2(n) = 3$ for $n\ge 3$.
\end{proposition}
\begin{proof}
That $R_2(n)\geq 3$ follows since any metric space on 3
points embeds isometrically in $\ell^2_2$. To show that $R_2(n) \le
3$, we construct a metric space on $n>3$
points, no $4$-point subspace of which
embeds isometrically in $\ell_2$. Fix an integer $n>3$ and let
$\{a_i\}_{i=0}^n$ be an increasing sequence such that $a_0=0$, $a_1=1$ and for $1\le
i<n$, $a_{i+1}\ge 2(n+1)a_i$. Fix some $0<\epsilon<1/(2a_n)$. It is
easily verified that
$d(i,j)=|i-j|-\epsilon a_{|i-j|}$ is a metric on $\{1,2,\ldots,n\}$. We
show that for
$\epsilon$ small enough no four points in $(\{1,\ldots,n\},d)$ embed
isometrically in $\ell_2$. Fix four integers $1\le i_1<i_2<i_3<i_4\le n$ and
set $j=i_2-i_1$, $k=i_3-i_2$, $l=i_4-i_3$. Suppose that for every $\epsilon>0$
there exists an isometric embedding $f:(\{i_1,i_2,i_3,i_4\},d) \rightarrow
\ell_2^3$. Without loss of generality we may assume that
$f(i_1)=(\alpha,\beta,\gamma)$, $f(i_2)=(0,0,0)$, $f(i_3)=(k- \epsilon
a_k,0,0)$ and $f(i_4)=(p,q,0)$. Then: {\setlength\arraycolsep{1pt}
\begin{eqnarray*}
2\alpha(k-\epsilon a_k)&=&2\langle
f(i_1),f(i_3)\rangle\\
&=&\|f(i_1)-f(i_2)\|_2^2+\|f(i_3)-f(i_2)\|_2^2-\|f(i_3)-f(i_1)\|_2^2\\
&=& (j-\epsilon a_j)^2+(k-\epsilon a_k)^2-(j+k-\epsilon
a_{j+k})^2.
\end{eqnarray*}}
Hence,
$$
\alpha\leq
-j+\frac{\epsilon}{k}[(k+j)a_{k+j}-ja_j-ka_k-ja_k]+O(\epsilon^2).
$$
Similarly:
$$
p \geq (k+l)+\frac{\epsilon}{k}[(k+l)a_k-(k+l)a_{k+l}-ka_k+la_l]+O(\epsilon^2).
$$
Now:
{\setlength\arraycolsep{1pt}
\begin{eqnarray*}
j+k+l&-&\epsilon a_{j+k+l}=
\|f(i_4)-f(i_1)\|_2\\
&\ge&
p-\alpha\\
&\geq&j+k+l+\frac{\epsilon}{k}[(k+l)a_k-(k+l)a_{k+l}+la_l-(k+j)a_{k+j}+ja_j+ja_k]+O(\epsilon^2).
\end{eqnarray*}}
Letting $\epsilon$ tend to zero we deduce that:
{\setlength\arraycolsep{1pt}
\begin{eqnarray*}
a_{j+k+l}\le
\left(1+\frac{j}{k}\right)a_{k+j}+\left(1+\frac{l}{k}\right)a_{k+l}-\frac{l}{k}a_l-\frac{j}{k}a_j-
\frac{j+k+l}{k}a_k<2(n+1)a_{j+k+l-1},
\end{eqnarray*}}
which is a contradiction.
\end{proof}

The argument above is quite specific to $\ell_2$, and so we now consider any
uniformly convex normed space. The modulus of
uniform convexity of a normed space $X$ is defined by:

$$
\delta_X(\epsilon)=\inf\left\{1-\frac{\|a+b\|}{2};\ \|a\|,\|b\|\le 1
\quad \mathrm{and}\quad \|a-b\|\ge \epsilon\right\}.
$$
$X$ is said to be uniformly convex if $\delta_X(\epsilon)>0$ for
every $0<\epsilon\le 2$. The $L_p$ spaces $1<p<\infty$, are known
to be uniformly convex. For a uniformly convex space $X$,
$\delta_X$ is known to be continuous and strictly increasing on
$(0,2]$.

Assume that $X$ is a uniformly convex normed space and $a,b\in X \setminus\{0\}$. Then:
{\setlength\arraycolsep{1pt}
\begin{eqnarray*}
\left\|\frac{a}{\|a\|}+\frac{b}{\|b\|}\right\| &=&
\left\|\left(\frac{1}{\|a\|}+\frac{1}{\|b\|}\right)(a+b)-\frac{a}{\|b\|}-\frac{b}{\|a\|}
\right\|\\
&\ge& \left(\frac{1}{\|a\|}+\frac{1}{\|b\|}\right)\|a+b\|-\frac{\|a\|}{\|b\|}-\frac{\|b\|}
{\|a\|}\\
&=&2-\left(\frac{1}{\|a\|}+\frac{1}{\|b\|}\right)(\|a\|+\|b\|-\|a+b\|).
\end{eqnarray*}}
Now,
{\setlength\arraycolsep{1pt}
\begin{eqnarray*}
\delta_X\left(\left\|\frac{a}{\|a\|}-\frac{b}{\|b\|}\right\|\right)
\le 1-\frac{1}{2}\cdot
\left\|\frac{a}{\|a\|}+\frac{b}{\|b\|}\right\|
\le\frac{1}{2}\cdot\left(\frac{1}{\|a\|}+\frac{1}{\|b\|}\right)(\|a\|+\|b\|-\|a+b\|).
\end{eqnarray*}}
Hence
$$
\left\|\frac{a}{\|a\|}-\frac{b}{\|b\|}\right\|\le\delta_X^{-1}\left(
\frac{1}{2}\cdot\left(\frac{1}{\|a\|}+\frac{1}{\|b\|}\right)(\|a\|+\|b\|-\|a+b\|)\right).
$$
Take $x,y,z\in X$ and apply this inequality for $a=x-y$, $b=y-z$. It follows that:
\begin{eqnarray}\label{eq:use}
\nonumber\left\|y-\left(\frac{\|y-z\|}{\|x-y\|+\|y-z\|}\right.\right.\!\!\!\!\!\!\!&\cdot&\!\!\!\!\!\!\!\left.\left.
x+\frac{\|x-y\|}{\|x-y\|+\|y-z\|}
\cdot z\right)\right\| \\
&\le& \frac{\|x-y\|\cdot\|y-z\|}{\|x-y\|+\|y-z\|}\cdot
\delta_X^{-1}
\left(\frac{\|x-y\|+\|y-z\|-\|x-z\|}{\min\{\|x-y\|,\|y-z\|\}}\right).
\end{eqnarray}
This inequality is the way uniform convexity is going to be applied in the sequel. Indeed, we
have the following ``metric'' consequence of it:
\begin{lemma}\label{convexity}
Let $X$ be a uniformly convex normed space and $x_1,x_2,x_3,x_4\in X$ be distinct. Then:
{\setlength\arraycolsep{1pt}
\begin{eqnarray*}
\frac{\|x_1-x_2\|+\|x_2-x_3\|-\|x_1-x_3\|}{2\|x_2-x_3\|}
&\le&\delta_X^{-1}
\left(\frac{\|x_1-x_3\|+\|x_3-x_4\|-\|x_1-x_4\|}{\min\{\|x_1-x_3\|,\|x_3-x_4\|\}}\right)
+\\
&\phantom{\le}&+ \delta_X^{-1}
\left(\frac{\|x_2-x_3\|+\|x_3-x_4\|-\|x_2-x_4\|}{\min\{\|x_2-x_3\|,\|x_3-x_4\|\}}\right).
\end{eqnarray*}}
\end{lemma}
Lemma~\ref{convexity} is a quantitative version of the fact that
in a uniformly convex space, if $\|x_1-x_4\|$ is approximately
$\|x_1-x_3\|+ \|x_3-x_4\|$ and $\|x_2-x_4\|$ is approximately
$\|x_2-x_3\|+ \|x_3-x_4\|$ then $\|x_1-x_3\|$ is approximately
$\|x_1-x_2\|+ \|x_2-x_3\|$. This fact is geometrically evident
since the first assumption implies that $x_3$ is almost on the
line segment connecting $x_1$ and $x_4$ and $x_2$ is almost on the
line segment connecting $x_1$ and $x_3$. It follows that $x_2$ is
almost on the line segment connecting $x_1$ and $x_3$, as
required. Since we are dealing with bi-Lipschitz embeddings, we
must formulate this phenomenon without referring to``line
segments".

\begin{proof}[Proof of Lemma~\ref{convexity}]
Define:
$$
\lambda=\frac{\|x_3-x_4\|}{\|x_1-x_3\|+\|x_3-x_4\|}\quad \mathrm{and}\quad
\mu=\frac{\|x_3-x_4\|}{\|x_2-x_3\|+\|x_3-x_4\|}.
$$
An application of \eqref{eq:use} twice gives:
$$
\|x_3-(\lambda x_1+(1-\lambda) x_4)\|\le \frac{\|x_1-\|x_3\|\cdot\|x_3-x_4\|}
{\|x_1-x_3\|+\|x_3-x_4\|}
\cdot \delta_X^{-1}
\left(\frac{\|x_1-x_3\|+\|x_3-x_4\|-\|x_1-x_4\|}{\min\{\|x_1-x_3\|,\|x_3-x_4\|\}}\right),
$$
and
$$
\|x_3-(\mu x_2+(1-\mu) x_4)\|\le
\frac{\|x-2-x_3\|\cdot\|x_3-x_4\|}{\|x_2-x_3\|+\|x_3-x_4\|} \cdot
\delta_X^{-1}
\left(\frac{\|x_2-x_3\|+\|x_3-x_4\|-\|x_2-x_4\|}{\min\{\|x_2-x_3\|,\|x_3-x_4\|\}}\right).
$$
By symmetry, we may assume without loss of generality that $\lambda\le \mu$. Now,
{\setlength\arraycolsep{1pt}
\begin{eqnarray*}
\left\|x_2-\frac{\lambda(1-\mu)}{\mu(1-\lambda)}x_1-
\frac{\mu-\lambda}{\mu(1-\lambda)}x_3\right\|
&=&\frac{1}{\mu}\left\|\mu
x_2+(1-\mu)x_4-x_3+\frac{1-\mu}{1-\lambda}(x_3-\lambda x_1-
(1-\lambda)x_4)\right\|\\
&\le&
\frac{1}{\mu}\|x_3-\mu x_2-(1-\mu)x_4\|+\frac{1-\mu}{\mu(1-\lambda)}\cdot\|x_3-\lambda x_1-
(1-\lambda)x_4\|\\
&\le&
\frac{\|x_2-x_3\|+\|x_3-x_4\|}{\|x_3-x_4\|}\cdot\frac{\|x_2-x_3\|\cdot\|x_3-x_4\|}
{\|x_2-x_3\|+\|x_3-x_4\|}\cdot \\
&\phantom{\le}&\cdot\delta_X^{-1}
\left(\frac{\|x_2-x_3\|+\|x_3-x_4\|-\|x_2-x_4\|}{\min\{\|x_2-x_3\|,\|x_3-x_4\|\}}\right)+\\
&\phantom{\le}&+\frac{\|x_2-x_3\|}{\|x_3-x_4\|}\frac{\|x_1-x_3\|+\|x_3-x_4\|}{\|x_1-x_3\|}
\frac{\|x_1-\|x_3\|\cdot\|x_3-x_4\|}
{\|x_1-x_3\|+\|x_3-x_4\|}
\cdot \\
&\phantom{\le}&\cdot\delta_X^{-1}
\left(\frac{\|x_1-x_3\|+\|x_3-x_4\|-\|x_1-x_4\|}{\min\{\|x_1-x_3\|,\|x_3-x_4\|\}}\right)\\
&=&
\|x_2-x_3\|\delta_X^{-1}
\left(\frac{\|x_1-x_3\|+\|x_3-x_4\|-\|x_1-x_4\|}{\min\{\|x_1-x_3\|,\|x_3-x_4\|\}}\right)+ \\
&\phantom{\le}&+\|x_2-x_3\|\delta_X^{-1}
\left(\frac{\|x_2-x_3\|+\|x_3-x_4\|-\|x_2-x_4\|}{\min\{\|x_2-x_3\|,\|x_3-x_4\|\}}\right)
.
\end{eqnarray*}}
Additionally,
{\setlength\arraycolsep{1pt}
\begin{eqnarray*}
\|x_2-x_1\|  &\le&
\left\|x_2-\frac{\lambda(1-\mu)}{\mu(1-\lambda)}x_1-\frac{\mu-\lambda}
{\mu(1-\lambda)}x_3\right\|+
\left\|x_1-\frac{\lambda(1-\mu)}{\mu(1-\lambda)}x_1-\frac{\mu-\lambda}{\mu(1-\lambda)}x_3
\right\|\\
 &=& \left\|x_2-\frac{\lambda(1-\mu)}{\mu(1-\lambda)}x_1-\frac{\mu-\lambda}
{\mu(1-\lambda)}x_3\right\| +\frac{\mu-\lambda}
{\mu(1-\lambda)}
\|x_1-x_3\|,
\end{eqnarray*}}
and
{\setlength\arraycolsep{1pt}
\begin{eqnarray*}
\|x_2-x_3\|  &\le&
\left\|x_2-\frac{\lambda(1-\mu)}{\mu(1-\lambda)}x_1-\frac{\mu-\lambda}
{\mu(1-\lambda)}x_3\right\|+
\left\|x_3-\frac{\lambda(1-\mu)}{\mu(1-\lambda)}x_1-\frac{\mu-\lambda}
{\mu(1-\lambda)}x_3\right\|\\
 &=&\left\|x_2-\frac{\lambda(1-\mu)}{\mu(1-\lambda)}x_1-\frac{\mu-\lambda}
{\mu(1-\lambda)}x_3\right\|+\frac{\lambda(1-\mu)}
{\mu(1-\lambda)}
\|x_1-x_3\|.
\end{eqnarray*}}
Summing up these estimates gives the required result.
\end{proof}

We can now prove the main result of this section:

\begin{theorem}\label{th:convex}
Let $X$ be a uniformly convex normed space with $dim(X)\ge 2$. Then for every
$n\ge 3$, $R_{X}(n)=3$. Moreover, for every $\delta:(0,2]\to(0,\infty)$ which
is continuous, increasing and $\delta\le \delta_{\ell_2}$, let $UC_{\delta}$
be the class of all normed spaces $X$ with $\delta_X\ge \delta$. Then for
each $n\ge 3$ there is a constant $\epsilon_n(\delta)>0$ such that $R_{
UC_{\delta}}({1+\epsilon_n(\delta)},n)=3$.
\end{theorem}

The proof of Theorem~\ref{th:convex} proceeds by constructing a
space in which each quadruple violates the conclusion of
Lemma~\ref{convexity}. The construction is done iteratively, by
adding one point at a time.

\begin{proof}[Proof of Theorem~\ref{th:convex}] That $R_{X}(n)\ge 3$ follows since
any $3$ point metric embeds isometrically into any $2$ dimensional normed space,
by a standard continuity argument.

Fix some $\delta:(0,2]\to (0,\infty)$ which is continuous, increasing and
$\delta\le \delta_{\ell_2}$. We shall construct inductively a sequence
$\{M_n\}_{n=3}^{\infty}$ of
metric spaces and numbers $\{\eta_n\}_{n=3}^{\infty}$ such that:

\medskip
\noindent{\bf a)} For every $n\ge 3$, $\eta_n>0$. Each $M_n$ is a metric on $\{1,\ldots,n\}$, and we denote $d^n_{ij}=
d_{M_n}(i,j)$.

\medskip
\noindent{\bf b)} For every $1\le i<j<k \le n$,
$$
d_{i,j}^n+d_{jk}^n-d_{i,k}^n-\eta_n\ge
2d_{j,k}^n\left[\delta^{-1}\left(\frac{d_{i,k}^n+d_{k,n}^n-
d_{i,n}^n}{\min\{d_{i,k}^n,d_{k,n}^n\}}\right)+\delta^{-1}\left(\frac{d_{j,k}^n+d_{k,n}^n-
d_{j,n}^n}{\min\{d_{j,k}^n,d_{k,n}^n\}}\right)\right].
$$
Lemma~\ref{convexity} immediately implies that there is
a constant $\epsilon_n(\delta)>0$ such that for every $1\le i<j<k<l\le n$ and
for every normed space $X$ with $\delta_X\ge \delta$:
$$
c_X(\{i,j,k,l\},d_{M_n})\ge 1+\epsilon_n(\delta),
$$
as required.

$M_3$ is the equilateral metric on $\{1,2,3\}$, in which case
$\eta_3=1$. We construct $M_{n+1}=(\{1,\ldots,n+1\},d^{n+1})$
as an extension of $M_n$, by setting
$$
d^{n+1}_{n,n+1}=1-s/2 \quad \text{and} \quad \forall 1\le i< n, \ d^{n+1}_{i,n+1}=d^n_{in}+
1-s.
$$
This is indeed a definition of a metric as long as
$0< s \leq \min \{1,2\min_{1\leq i<n}d^n_{i,n}\}$ (this
fact follows from a simple case analysis).

We are left to check condition ${\bf b)}$.
Fix $1\le i<j<k\le n$. If $k\neq n$ then:
{\setlength\arraycolsep{1pt}
\begin{eqnarray*}
d_{i,j}^{n+1}+d_{j,k}^{n+1}-d_{i,k}^{n+1}-
\eta_n&=&d_{i,j}^n+d_{j,k}^n-d_{i,k}^n-\eta_n\\
&\ge&
2d_{j,k}^n\left[\delta^{-1}\left(\frac{d_{i,k}^n+d_{k,n}^n-
d_{i,n}^n}{\min\{d_{i,k}^n,d_{k,n}^n\}}\right)+\delta^{-1}\left(\frac{d_{j,k}^n+d_{k,n}^n-
d_{j,n}^n}{\min\{d_{j,k}^n,d_{k,n}^n\}}\right)\right]\\
&\ge&
2d_{j,k}^n\left[\delta^{-1}\left(\frac{d_{i,k}^n+(d_{k,n}^n+1-s)-
(d_{i,n}^n+1-s)}{\min\{d_{i,k}^n,d_{k,n}^n+1-s\}}\right)\right.\\
&\phantom{\le}&\left.+\delta^{-1}
\left(\frac{d_{j,k}^n+(d_{k,n}^n+1-s)-
(d_{j,n}^n+1-s)}{\min\{d_{j,k}^n,d_{k,n}^n+1-s\}}\right)\right]\\
&=&
2d_{j,k}^{n+1}\left[\delta^{-1}\left(\frac{d_{i,k}^{n+1}+d_{k,n+1}^{n+1}-
d_{i,n+1}^{n+1}}{\min\{d_{i,k}^{n+1},d_{k,n+1}^{n+1}\}}\right)+
\delta^{-1}\left(\frac{d_{j,k}^{n+1}+d_{k,n+1}^{n+1}-
d_{j,n+1}^{n+1}}{\min\{d_{j,k}^{n+1},d_{k,n+1}^{n+1}\}}\right)\right].
\end{eqnarray*}}

It remains to check ${\bf b)}$
for the quadruple $\{i,j,n,n+1\}$. Condition ${\bf b)}$ for $M_n$ implies that:
$$
d_{ij}^{n+1}+ d_{jn}^{n+1}-d_{in}^{n+1}\ge \eta_n.
$$
On the other hand,
\begin{multline*}
2d_{j,n}^{n+1}\left[\delta^{-1}\left(\frac{d_{i,n}^{n+1}+d_{n,n+1}^{n+1}-
d_{i,n+1}^{n+1}}{\min\{d_{i,n}^{n+1},d_{n,n+1}^{n+1}\}}\right)+\delta^{-1}
\left(\frac{d_{j,n}^{n+1}+d_{n,n+1}^{n+1}-
d_{j,n+1}^{n+1}}{\min\{d_{j,n}^{n+1},d_{n,n+1}^{n+1}\}}\right)\right]=\\
=2d_{j,n}^n\left[\delta^{-1}\left(\frac{s/2}{\min\{d_{i,n}^n,1-s/2\}}\right)+\delta^{-1}
\left(\frac{s/2}{\min\{d_{j,n}^n,1-s/2\}}\right)\right],
\end{multline*}
so that condition ${\bf b)}$ will hold
when $s$ is small enough such that the quantity
above is at most $\eta_n/2$ and with $\eta_{n+1}=\eta_n/2$.
\end{proof}

\begin{corollary}\label{cor:isometric-lp}
For all $1<p<\infty$, $R_p(n) = 3$ for $n\ge 3$.
\end{corollary}
\bigskip

We end this section with a simple lower bound for the isometric
Ramsey problem for graphs. We do not know the
asymptotically tight bound in this setting.


\begin{proposition}
Let $G$ be an unweighted graph of order $n$. Then there is a set
of $\Omega\left(\sqrt{\frac{\log n}{\log \log n}}\right)$ vertices
in $G$ whose metric embeds isometrically into $\ell_2$.
\end{proposition}
\begin{proof}
Let $\Delta$ be the diameter of $G$. The
shortest path between two diameterically far vertices is isometrically embeddable in
$\ell_2$. On the other hand, the Bourgain,Figiel,
Milman theorem \cite{bfm}
yields that for every $0<\epsilon<1$ a subset $N\subset V$ which is
$(1+\epsilon)$ embeddable in Hilbert space and
$|N|=\Omega\left(\frac{\epsilon}{\log(2/\epsilon)}\log n \right)$.
When $\epsilon=\frac{1}{2\Delta}$, such an embedding is an isometry.
Hence we can
always extract a subset of $V$ which is isometrically embeddable
in $\ell_2$ with cardinality
$$
\Omega\left(\max\left\{\Delta,\frac{\log n}{\Delta\log \Delta}\right\}\right)=\Omega
\left(\sqrt{\frac{\log n}{\log \log n}}\right),
$$
as claimed.
\end{proof}

\bigskip
\noindent{\bf Acknowledgments:} The authors would like to express
their gratitude to Guy Kindler for some helpful discussions.

\bibliographystyle{plain}

\bigskip
\bigskip

\noindent Yair Bartal, Institute of Computer Science, Hebrew
University, Jerusalem 91904, Israel. \\{\bf yair@cs.huji.ac.il}

\medskip
\noindent Nathan Linial, Institute of Computer Science, Hebrew
University, Jerusalem 91904, Israel. \\ {\bf nati@cs.huji.ac.il}

\medskip
\noindent Manor Mendel, Institute of Computer Science, Hebrew
University, Jerusalem 91904, Israel. \\ {\bf
mendelma@cs.huji.ac.il}

\medskip
\noindent Assaf Naor, Theory Group, Microsoft Research, One
Microsoft Way 113/2131, Redmond WA 98052-6399, USA. \\ {\bf
anaor@microsoft.com}

\bigskip
\bigskip
\noindent 2000 AMS Mathematics Subject Classification: 52C45,
05C55, 54E40, 05C12, 54E40.

\end{document}